\documentclass{amsart}

\usepackage{amscd}
\usepackage{amsfonts}
\usepackage{amssymb}
\usepackage{latexsym}

\newcommand{\ncm}{\newcommand}

\ncm{\beq}{\begin{equation}}
\ncm{\eeq}{\end{equation}}
\ncm{\bea}{\begin{eqnarray}}
\ncm{\eea}{\end{eqnarray}}
\ncm{\beanon}{\begin{eqnarray*}}
\ncm{\eeanon}{\end{eqnarray*}}

\newtheorem{thm}{Theorem}[section]
\newtheorem{pro}[thm]{Proposition}
\newtheorem{lem}[thm]{Lemma}

\newtheorem{lem&def}[thm]{Lemma \& Definition}

\theoremstyle{definition}
\newtheorem{defi}[thm]{Definition}
\newtheorem{exa}[thm]{Example}

\theoremstyle{remark}

\numberwithin{equation}{section}

\def\Vec{\mathsf{Vec}}

\def\Alg{\mathsf{Alg}}
\def\Cat{\mathsf{Cat}}

\ncm{\rep}{\mathsf{rep}}
\ncm{\amp}{\mathsf{amp}}
\def\C{\mathcal{C}}
\def\M{\mathcal{M}}

\ncm{\End}{\operatorname{End}}

\def\Hom{\mbox{\rm Hom}\,}

\def\id{\mbox{\rm id}\,}
\def\Center{\mbox{\rm Center}\,}
\ncm{\Ind}{\operatorname{Ind}}

\def\o{\otimes}    
\ncm{\oo}{\Box\,}  
\ncm{\oh}{\x}      
\ncm{\ov}{\circ}   
\def\x{\times}     
\def\bra{\langle}
\def\ket{\rangle}
\def\under{\mbox{\rm\_}\,}
\ncm{\rarr}[1]{\stackrel{#1}{\longrightarrow}}
\ncm{\larr}[1]{\stackrel{#1}{\longleftarrow}}
\def\cop{\Delta}
\def\ducop{\hat\cop}

\def\eps{\varepsilon}
\def\dueps{\hat\eps}
\def\du1{\hat 1}
\def\1{_{(1)}}
\def\2{_{(2)}}
\def\3{_{(3)}}
\def\tr{\mbox{\rm tr}}
\def\PL{\pi^{\scriptscriptstyle L}}
\def\PR{\pi^{\scriptscriptstyle R}}
\def\duA{\hat A}
\def\la{\!\rightharpoonup\!}
\def\ra{\!\leftharpoonup\!}
\ncm{\ma}[1]{#1\triangleright}
\def\du1{\hat 1}
\def\F{\mathcal{F}}

\begin{document}

\title[Weak Hopf algebras]{Weak Hopf algebra symmetries of $C^*$-algebra
inclusions}   
\author[K. Szlach\'anyi]{Korn\'el Szlach\'anyi} 
\address{Research Institute for Particle and Nuclear Physics, Budapest\\
 H-1525 Budapest, P. O. Box 49, Hungary}
\email{szlach@rmki.kfki.hu} 
\thanks{This paper is based on a series of talks given by the author on the
Conference "Index Theory and Physics" at the Universit\'a Degli Studi di
Bologna in February, 2000. The paper partly contains material that has not
been published yet. 
} 
\maketitle

In 2-dimensional conformal field theory or in 1-dimensional quantum lattice
systems the superselection sectors (charges) $p,q,\dots$ often have
composition laws $p\o q=\sum_r\,N_{pq}^r\cdot r$ such that
the non-negative integer multiplicities $N_{pq}^r$ do not admit an
integer solution for the dimensions $d_r$ of the equation
$\sum_r\,N^r_{pq}\,d_r\,=\,d_pd_q$. This implies that it is impossible to find
a semisimple quantumgroup, more precisely a semisimple Hopf algebra, $H$ the
irreducible representations of which would be in one-to-one correspondence
with the sectors $p,q,\dots$ and the tensor product of irreducible
representations of $H$ would produce precisely the numbers $N_{pq}^r$ as their
branching numbers. 

A solution for this problem has been found in \cite{BSz} where we proposed to
replace $C^*$-Hopf algebras with $C^*$-{\em weak Hopf algebras} the
comultiplication of which is coassociative but not unit preserving. This
allows the dimension of a tensor product of representations be smaller than
the product of the dimensions of these representations,
\[
\dim(D\o D')\ \leq \dim D\cdot \dim D'\ .
\]
$C^*$-weak Hopf algebras extended earlier generalizations by T. Yamanouchi
\cite{Yam} and Y. Hayashi \cite{Hay} and can be considered as an axiomatic
approach to the depth 2 paragoups of A. Ocneanu. A detailed analysis of the
structure of weak Hopf algebras has been carried out in \cite{BNSz, Nill, NV1,
BSz2}. For a recent review see \cite{NV4}. 

Independently of the above development the motivation to axiomatize a fairly
non-commutative Poisson groupoid has lead J-H. Lu to the introduction of
{\em Hopf algebroids} in \cite{Lu}. Later it was shown by P. Etingof
and D. Nikshych \cite{EN} that weak Hopf algebras are just the finite
dimensional versions of Hopf algebroids. See also \cite{Szlach}. 

The definition of a Hopf algebroid $A$ contains the data of a base algebra $B$
and two algebra maps $s\colon B^{op}\to A$ and $t\colon B\to A$, called the
source and the target, respectively. They are non-commutative versions of the
algebra of functions over the space of units of a groupoid. In a weak Hopf
algebra $A$ the images of these maps has been called $A^R$ and $A^L$,
respectively, although they are not part of the data but have been
"discovered" in studying the properties of the coproduct $\cop\colon A\to
A\o A$.   

$C^*$-weak Hopf algebras can be applied to the characterization of finite
index depth 2 inclusions $N\subset M$ of von Neumann algebras. If the centers
of $M$ and $N$ are finite dimensional there exists a $C^*$-weak Hopf algebra
$A$ acting regularly on $M$ such that $N$ is the invariant subalgebra. Such an
approach has been initiated in \cite{NSzW} and worked out for type II$_1$
factors by D. Nikshych and L. Vainerman in \cite{NV2, NV3}. One should mention
here an other approach \cite{EV, E2} by M. Enock and J.-M. Vallin which uses
Hopf bimodules \cite{Val1} to characterize depth 2 inclusions with arbitrary
(non-finite) index.

The paper is organized as follows. In Section \ref{s WHA} the axioms of weak
bialgebras  and weak Hopf algebras are presented with a brief description of
the properties of $A^L$, $A^R$, representation theory, the $C^*$-structure,
and the Haar measure. In Section \ref{s modalg} we summarize the requirements
for a "regular" action of a $C^*$-weak Hopf algebra on a unital $C^*$-algebra.
The first two sections are based on earlier papers on the subject while Section
\ref{s reco} contains new results. In Section \ref{s reco} we outline the proof
of a reconstruction theorem stating that if a certain inclusion $N\subset M$ of
$C^*$-algebras is given then there is a $C^*$-weak Hopf algebra $A$ and a
regular action of $A$ on $M$ such that $N$ is the invariant subalgebra. We
actually work with abstract 2-categorical "inclusions".

\section{Weak Hopf algebras} \label{s WHA}
\subsection{Weak bialgebras}
\begin{defi}
A weak bialgebra over a field $K$ consists of the data $\bra
A,m,u,\cop,\eps\ket$ where
\begin{enumerate}
\item $\bra A,m,u\ket$ is an algebra, i.e., the multiplication $m\colon A\o
A\to A$ and the counit $u\colon K\to A$ satisfy
\begin{enumerate} 
\item associativity: $m\circ(m\o \id )=m\circ(\id \o m)$ and
\item unit properties: $m\circ(u\o\id )=\id =m\circ(\id \o u)$
\end{enumerate}
\item $\bra A,\cop,\eps\ket$ is a coalgebra, i.e., the
comultiplication $\cop\colon A\to A\o A$ and the counit $\eps\colon A\to K$
satisfy
\begin{enumerate}
\item coassociativity: $(\cop\o\id )\circ\cop=(\id \o\cop)\circ\cop$ and
\item counit properties: $(\eps\o\id )\circ\cop=\id =(\id \o\eps)\circ\cop$
\end{enumerate}
\item the algebra and coalgebra structures satisfy the compatibility conditions
\begin{enumerate}
\item $\cop$ is multiplicative / $m$ is comultiplicative, i.e., as maps $A\o
A\to A\o A$, 
$$\cop\circ m=(m\o m)\circ(\id \o\tau\o\id )\circ(\cop\o\cop)$$
where $\tau\colon A\o A\to A\o A$ denotes the flip map $x\o y\mapsto y\o x$,
\item $\eps$ is weakly multiplicative, i.e., as maps $A\o A\o A\to K$,
\beanon
(\eps\o \eps)\circ(m\o m)\circ(\id \o\cop\o\id )&=&\eps\circ
m\circ(m\o\id )\\ 
(\eps\o \eps)\circ(m\o m)\circ(\id \o\cop^{op}\o\id )&=&\eps\circ
m\circ(m\o\id ) 
\eeanon
where $\cop^{op}:=\tau\circ\cop$ is opposite comultiplication,
\item $u$ is weakly comultiplicative, i.e., as maps $K\to A\o A\o A$,
\beanon
(\id \o m\o\id )\circ(\cop\o\cop)\circ(u\o u)&=&(\cop\o\id )\circ\cop\circ
u\\
(\id \o m^{op}\o\id )\circ(\cop\o\cop)\circ(u\o
u)&=&(\cop\o\id )\circ\cop\circ u
\eeanon
where $m^{op}:=m\circ\tau$ is opposite multiplication.
\end{enumerate}
\end{enumerate}
\end{defi}

Usually one uses elements of $A$ to express identities and this is what we
will do in most of the cases. So we shall write $xy$ instead of $m(x\o y)$ and
use the unit element $1$ instead of $u\colon k\mapsto 1k$. Then axiom (3.a)
could then be simply written as $\cop(ab)=\cop(a)\cop(b)$. Or, using Sweedler's
notation $\cop(a)=a\1\o a\2$ which suppresses a possible summation, (3.a)
takes the form $(ab)\1\o(ab)\2=a\1 b\1\o a\2 b\2$, well, the beauty of which is
not very convincing. However, the weak multiplicativity properties under (3.b)
can be rather nicely written using Sweedler's notation as
$\eps(ab\1)\eps(b\2 c)=\eps(abc)=\eps(ab\2)\eps(b\1 c)$.  

Anyhow, the above set of axioms show manifestly the selfduality
of the structure. Considering the axioms as commutative diagrams in the
category $\Vec K$ of vector spaces then the axioms stay invariant if we
reverse the directions of all the arrows and change simultaneously
$\cop \leftrightarrow m$ and $\eps \leftrightarrow u$.

If $A$ is a finite dimensional bialgebra then $\hat A:=\Hom(A,K)$ is also a
bialgebra if we define the structure maps $\hat m,\ducop,\hat u,\dueps$ by
means of the canonical pairing $\bra\under,\under\ket\colon\hat A\x A\to K$, 
\bea
\bra \varphi\psi,a\ket&:=&\bra \varphi\o \psi,\cop(a)\ket\\
\bra\ducop(\varphi),a\o b\ket&:=&\bra \varphi,ab\ket\\
\bra\hat 1,a\ket&:=&\eps(a)\\
\dueps(\varphi)&:=&\bra \varphi,1\ket\,,
\eea
where $\varphi,\psi\in\duA$, $a,b\in A$, and we used the notation
$\bra\under,\under\ket$ also for the pairing of $A\o A$ with
$\duA\o\duA\cong\Hom(A\o A,K)$. This is exactly how ordinary bialgebras behave.
A weak bialgebra $A$ becomes a bialgebra in case either of the following
conditions hold:
\begin{enumerate}
\item $\cop$ is unital, $\cop(1)=1\o 1$, 
\item $\eps$ is multiplicative, $\eps(ab)=\eps(a)\eps(b)$ for all $a,b\in A$.
\end{enumerate}

The basic feature in weak bialgebras (WBA's for short) is the presence of
two canonical subalgebras $A^L\subset A$ and $A^R\subset A$ both
of which reduce to $K$ in case of $A$ is an ordinary bialgebra. The following
are equivalent definitions of $A^L$ and $A^R$.
\begin{lem&def}
The following conditions on elements $l$, respectively $r$, of a WBA $A$ are
equivalent
\begin{enumerate}
\item $\cop(l)=(l\o 1)\cop(1)=\cop(1)(l\o 1)$
      \newline $\cop(r)=(1\o r)\cop(1)=\cop(1)(1\o r)$
\item $(\eps\o\id)(\cop(1)(l\o 1))=l$
      \newline $(\id\o\eps)((1\o r)\cop(1))=r$
\item $\exists f\colon A\to K$ such that $(f\o \id)\circ\cop(1)=l$\label{leg}
      \newline $\exists f\colon A\to K$ such that $(\id\o f)\circ\cop(1)=r$
\end{enumerate}
The set of such $l$'s and $r$'s will be denoted by $A^L$ and $A^R$,
respectively.
\end{lem&def}

It not difficult to see from the first of these three properties that $A^L$ and
$A^R$ are subalgebras of $A$ and that they commute with each other. $\cop(1)$
is an idempotent which not only belongs to $A^R\o A^L$ by property (2)
but also spans $A^L$ and $A^R$ in the sense of property (3).  
For the proof of these properties and many more we refer to \cite{BNSz, Nill}.

\subsection{Representations of WBA's} \label{ss: rep WBA}
The meaning of $A^L$, $A^R$ can be better understood in terms of the
representation theory of the weak bialgebra. Let $_A\M$ denote the category of
left $A$-modules. We define a monoidal structure on this category using the
coalgebra structure of $A$. The monoidal product $V\oo W$ of two $A$-modules
is defined as the subspace of the $K$-module tensor product $V\o W$
which is the image of the idempotent $\cop(1)$. The monoidal product of
intertwiners $t\colon V\to V'$ and $s\colon W\to W'$ is the restriction of
$t\o s\colon V\o W\to V'\o W'$ to the subspace $V\oo W=1\1\cdot V\o 1\2\cdot
W$.

This procedure is similar to how representation theory of a bialgebra or Hopf
algebra is related to the representation theory of the underlying field or
ring $K$. The "only" difference is the need of projecting out a subspace by
acting with $\cop(1)$. There is another, more elegant way of formulating the
monoidal structure of $_A\M$ which uses the fact that every bimodule category
has a very natural monoidal structure. 

Let $L$ denote the ring $A^L$ and its injection to $A$ be $t\colon
L\to A$. The ring $A^R$ is isomorphic to $L$ with opposite multiplication
via the map
\beq
s\colon L^{op}\to A\,,\quad s(l)=1\1\eps(1\2 l)
\eeq
Since the images of $s$ and $t$ commute, we can make $A$ into an
$L$-$L$-bimodule with the definitions
\beq
l\in L,\ a\in A\ \mapsto\ l\cdot a:=t(l)a\,,\quad a\cdot l:=s(l)a\,.
\eeq
Then we have the forgetful functor $\phi\colon\,_A\M\to\,_L\M_L$ from the
category of $A$-modules to the {\em monoidal} category of $L$-$L$-bimodules.
Thus a monoidal structure on $_A\M$ can be introduced by requireing that the
forgetful functor $\phi$ be strictly monoidal. 

This means precisely that we {\em define} the underlying space of the
$A$-module $V\oo W$ as the bimodule tensor product $V\o_L W$, which is indeed
smaller then the $K$-tensor product $V\o W$ and can be shown to be isomorphic
to the subspace projected out by $\cop(1)$, as above.
The proof uses the following identity shown in \cite{BNSz}, (2.31a),
\beq
a\1 s(l)\o a\2\ =\ a\1\o a\2 t(l)\,,\qquad l\in L,\ a\in A\,.
\eeq

Now it is not surprising that the monoidal unit $U$ of $_A\M$, also called the
{\em trivial representation} is not one-dimensional but is a representation on
the space $L=A^L$. So we have a left action of $A$ on $L$, but also, we had a
right action of $L$ on $A$, so the notation $a\cdot l$ would be ambiguous. 
Instead we use the special notation
\beq \label{eq: triv}
a\triangleright l\ :=\ \eps(1\1 al)1\2.
\eeq
In this way $L$ becomes a left $A$-module denoted $U$ for which there exist
natural isomorphisms $U\oo V\cong V\cong V\oo U$ for all $A$-modules $V$. Of
course, they are the natural isomorphisms $L\o_L V\cong V\cong V\o_L L$ of the
bimodule category $_L\M_L$. The very fact that these natural isomorphisms (and
also the associativity natural isomorphism) for the two categories $_A\M$ and
$_L\M_L$ coincide, means that the forgetful functor $_A\M\to\,_L\M_L$ is
strictly monoidal. 

The description of weak bialgebras through the bimodule category over $L$
yields a direct connection to the {\em bialgebroids} defined by Lu in
\cite{Lu}. A short review of the connection of these concepts can be found in
\cite{Szlach}.

An important example for a left $A$-module is the linear space of the dual
bialgebra $\duA$ endowed with the left action $a\la \varphi:=\varphi\1 
\bra \varphi\2,a\ket$,
where $a\in A$ and $\varphi\in\duA$. Similarly, a right $A$-module structure on
$\duA$ can be defined by the right action $\varphi\ra a:=\bra \varphi\1,a\ket
\varphi\2$. For showing that these are indeed actions one can use the
identities 
\beq
\bra b\la\varphi,a\ket\ =\ \bra\varphi,ab\ket\ =\ \bra\varphi\ra a,b\ket\,.
\eeq 

\subsection{Weak Hopf algebras}

A weak Hopf algebra (WHA for short) is a weak bialgebra $A$ for which there
exists an antipode, i.e., a linear map $S\colon A\to A$ satisfying the
antipode axioms
\bea
a\1 S(a\2)&=&\eps(1\1 a)1\2    \label{eq left antipode ax}\\
S(a\1)a\2&=&1\1\eps(a1\2)      \label{eq right antipode ax}\\
S(a\1)a\2 S(a\3)&=& a
\eea
for all $a\in A$. Denoting by $\PL$ and $\PR$ the right hand sides of axioms
(\ref{eq left antipode ax}) and (\ref{eq right antipode ax}), respectively, we
can make the following observations. In any WBA $\PL$ and $\PR$ are idempotents
in $\End A$ that project onto $A^L$ and $A^R$, respectively. But they are
idempotents also in the convolution sense. The convolution product is the
associative binary operation on $\End A$ defined by
\beq
f\ast g\ :=\ m\circ(f\o g)\circ\cop
\eeq
which has $I=u\circ\eps$ as its unit. One obtains the following identities in
any WBA:
\bea
\PL\ast\id_A=\id_A\,,&&\quad\PL\ast\PL=\PL\,,\label{eq target of id}\\
\id_A\ast\PR=\id_A\,,&&\quad \PR\ast\PR=\PR\,.\label{eq source of id}
\eea
In Hopf algebra theory one defines the
antipode as the convolution inverse of $\id_A$. Assume that
$\sigma\ast\id_A=I$. Then
$$
I=\sigma\ast\id_A=\sigma\ast\id_A\ast\PR=I\ast\PR=\PR
$$
but the image of $\PR$ is $A^R$ while the image of $I$ is $K$. So if the WBA
is not a bialgebra then $\id_A$ is never convolution invertible. Equations  
(\ref{eq target of id}) and (\ref{eq source of id}) suggest that a convolution
inverse can exist in a groupoid sense because $\PL$ is the target and $\PR$ is
the source projection of $\id_A$. Therefore we require for a convolution
inverse $S$ that
\beq
\id_A\ast S=\PL\,,\quad S\ast\id_A=\PR
\eeq
which are precisely the first two antipode axioms. The third one 
is then equivalent to either one of the equations
\beq
S\ast\PL=S\,,\qquad \PR\ast S=S
\eeq
that formulate the requirement that the source and target of $S$ be the target
and source of $\id_A$, respectively.

\begin{pro}[see Thm 2.10 of \cite{BNSz}]
The antipode $S$ of a WHA $A$ satisfies the following properties.
\begin{enumerate}
\item $S(ab)=S(b)S(a)$ for $a,b\in A$, i.e., $S$ is antimultiplicative,
\item $S(a)\1\o S(a)\2=S(a\2)\o S(a\1)$ for $a\in A$, i.e., $S$ is
anticomultiplicative,
\item $S(1)=1$, $\eps\circ S=\eps$,
\item $S(A^L)=A^R$ and $S(A^R)=A^L$,
\item if $\dim A <\infty$ then $S$ is invertible.
\end{enumerate}
\end{pro}

\begin{exa}
Let $A=KG$ be the groupoid algebra of a finite groupoid $G$. Define
comultiplication $\cop$ as the linear extension of the diagonal map $g\mapsto
g\o g$, $g\in G$. This map is not unital 
\beq
\cop(1)=\sum_{u\in O}\ u\o u \neq \sum_{u\in O}\sum_{v\in O} u\o v=1\o 1
\eeq
where $O\subset G$ denotes the set of units.
One can easily check that $\cop$ makes $KG$ into a WHA with
\beq
\eps(g)=1\,,\qquad S(g)=g^{-1}\,,\quad g\in G\,.
\eeq
In this WHA $A^L=A^R$ and they coincide with the groupoid algebra $KO$ of the
(totally disconnected) subgroupoid of units of $G$. Thus $A^L$ is Abelian
and spanned by the pairwise orthogonal idempotents $u\in O$.
\end{exa} 

\subsection{Weak $C^*$-Hopf algebras}

\begin{defi}
A $C^*$-WHA is a weak Hopf algebra $\bra A,\cop,\eps,S\ket$ over 
the complex field $\mathbb{C}$ in which $A$ is a finite
dimensional $C^*$-algebra and $\cop$ is a $^*$-algebra map.
\end{defi}
Since the counit and the antipode are uniquely determined by $\cop$, one can
immediately see that 
\beq
\overline{\eps(a)}=\eps(a^*)\,,\quad S(a^*)^*=S^{-1}(a)\,,\qquad a\in A\,.
\eeq
The identity $1\1^*\o 1\2^*\equiv \cop(1)^*=\cop(1^*)\equiv \cop(1)$ together
with (\ref{leg}) implies that $A^L$ and $A^R$ are closed under the
$^*$-operation, hence they are $C^*$-subalgebras of $A$. 

The dual WHA $\duA$ can be given a $^*$-operation by setting
\beq
\bra\varphi^*,a\ket\ :=\ \overline{\bra \varphi,S(a)^*\ket}\,.
\eeq
In this way $\duA$ becomes a $^*$-algebra and $\ducop$ a $^*$-homomorphism.
However, it is not obvious whether $\duA$ is also a $C^*$-algebra, inspite of
its finite dimension. What is missing is to show that $\duA$ possesses a
faithful $^*$-representation. What helps here is the existence of a Haar
measure. 
\begin{thm}
Let $A$ be a $C^*$-WHA. Then there exists a unique $h\in A$, called the Haar
integral (or Haar measure), such that 
\bea
ah=\PL(a)h\,,&\quad& ha=h\PR(a)\quad \forall a\in A\,,\\
\PL(h)\ =&1&=\ \PR(h)\,.
\eea
The Haar integral satisfies
\bea
h\1\o ah\2&=&S(a)h\1\o h\2\,\quad\forall a\in A\\
h\1 a\o h\2&=&h\1\o h\2 S(a)\,\quad\forall a\in A\\
h\la \duA&=&\duA^L\\
\duA\ra h&=&\duA^R\\
h\ =\ h^2&=&h^*\ =\ S(h)\\
\bra \varphi(\psi\ra a),h\ket&=&\bra(\varphi\ra S(a))\psi,h\ket\,,\quad\forall
     \label{rightinv}         \varphi,\psi\in\duA, \,\quad\forall a\in A\\
\bra (a\la\varphi)\psi,h\ket&=&\bra\varphi(S(a)\la\psi),h\ket\,,\quad\forall
     \label{leftinv}             \varphi,\psi\in\duA, \,\quad\forall a\in A
\eea
Moreover, the sesquilinear form on $\duA$ defined by
\beq
(\varphi,\psi)\ :=\ \bra\varphi^*\psi,h\ket
\eeq
is non-degenerate and positive, hence defines a Hilbert space structure on
$\duA$. It makes the left regular representation of $\duA$ into a faithful
$^*$-representation, in this way proving that $\duA$ is a $C^*$-WHA, too.
\end{thm}

Since $\dim A$ is finite, the Haar measure above is analogous to the Haar
measure on a finite group, i.e., the counting measure normalized to give a
total mass 1. 

Interpreting the pairing $\bra\psi, h\ket$ as the integral $\int\psi$ of the
function $\psi$ w.r.t. the Haar measure and the Sweedler arrows $a\la\under$
and $\under\ra a$ as left and right translations, respectively, properties
(\ref{rightinv}, \ref{leftinv}) become the right and left invariance of the
Haar measure, respectively. 

The maps $\varphi\mapsto (h\la \varphi)$ and $\varphi\mapsto(\varphi\ra h)$ are
conditional expectations from $\duA$ onto $\duA^L$ and $\duA^R$, respectively.
They tell us that $\duA^L$ is the space of functions invariant under left
translations by all $a\in A$. This space is, in general, different from the
space $\duA^R$ that consists of the functions invariant under right
translations.

The Haar state $\bra\under,h\ket$ on $\duA$ is not a trace, in general. In
order to compute its modular automorphism we have to acquaint with the
canonical grouplike element.

\begin{pro}
Let $A$ be a $C^*$-WHA and $h$ its Haar measure and let $\hat h$ be the Haar
measure of the dual WHA $\duA$. Then $\hat h\la h$ and $h\ra\hat h$ are
positive and invertible and we can define 
\beq
g\ :=\ g_Lg_R^{-1}\,\qquad \mbox{where}\ g_L=(\hat h\la h)^{1/2}\,,\quad
                     g_R=(h\ra\hat h)^{1/2}\,,
\eeq
and call it the canonical grouplike element of $A$.
It is the unique $g\in A$ for which
\bea
g&\geq &0\ \mbox{and invertible}\\
gag^{-1}&=&S^2(a)\,,\quad a\in A\\
\tr_r g&=&\tr_r g^{-1}\ \mbox{for all irrep $r$ of $A$}\,.
\eea
\end{pro}

The name "grouplike" refers to its invertibility and to the special form of
its coproduct 
\beq
\cop(g)\ =\ (g\o g)\cop(1)=\cop(1)(g\o g)\,.
\eeq
In Hopf algebras grouplike elements can be very sparse. The fact that $g$
exists in any $C^*$-WHA is related to that $g$ belongs to the "trivial"
subalgebra $A^LA^R$, which consists of the scalars in case of Hopf algebras.
Notice also that $g$ is not as grouplike as it could be, namely $g$ is not
unitary but positive instead.

The modular automorphism of the Haar state on $A$ can now be expressed in
terms of $g_L$ and $g_R$ as
\beq
\bra\hat h,ab\ket\ =\ \bra\hat h, b(g_Lg_R)a(g_Lg_R)^{-1}\ket\,,\quad a,b\in A
\,.
\eeq

One can prove that $\hat h$ is a trace precisely if $S^2=\id_A$. If this
happens the WHA is called a weak Kac algebra.

\subsection{Soliton sectors of $C^*$-WHA's}
The category $\rep A$ of finite dimensional $^*$-representations of a
$C^*$-WHA $A$ has been studied in detail in \cite{BSz2}. Here I would like to
emphasize only one aspect, the groupoid like vacuum structure of $\rep A$.
This property is independent of the fact that WHA's are quantum groupoids, it
reflects rather a "quantum 2-groupoid" feature. 

The monoidal unit of the monoidal category $\rep A$ is the GNS representation
associated to the positive linear functional $\eps\colon A\to \mathbb{C}$.
This is the representation mentioned in Subsection \ref{ss: rep WBA} endowed
with the scalar product $(l_1,l_2)=\eps(l_1^*l_2)$, $l_1,l_2\in A^L$ which
makes the left action (\ref{eq: triv}) a $^*$-representation of $A$. The point
is that this "trivial representation" may be reducible. This happens precisely
when the inclusion $A^L\subset A$ is not connected. The reason is that
$\End U\cong Z^L\equiv A^L\cap\Center A$ according to Proposition 2.15 of
\cite{BNSz}. The irreducibles $V_\mu$ occuring in the decomposition
\beq
U\ =\ \bigoplus_\mu\ V_\mu
\eeq
are called vacuum representations, the name borrowed from the theory of
solitons in 1+1-dimensional quantum field theory. Each of them is
selfconjugate and occurs with multiplicity 1 in $U$. WHA's with irreducible
trivial representation are called {\em pure} (since $\eps$ is pure) or
{\em connected} (since the Bratteli diagram of $A^L\subset A$ is connected).

If $V_r$ is an arbitrary
irreducible representation from the equivalence class (or sector, for short)
$r$ then there exists one and only one vacuum sector, called $r^L$, for
which  
\beq 
V_r\cong U\o V_r\cong\oplus_\mu V_\mu\o V_r\cong V_{r^L}\o V_r\,.
\eeq
Similarly, there exists one and only one vacuum sector $r^R$ for which $V_r\o
V_{r^R}\cong V_r$. The vacua $r^R$ and $r^L$ of the sector $r$ behave very much
like the source and target of a groupoid. 
\begin{itemize}
\item The monoidal product of sectors $p\o q=\{0\}$ if $p^R\neq q^L$.
\item If $p^R= q^L$ and $r$ is a sector contained in $p\o q$ then $r^L=p^L$
and $r^R=q^R$. 
\item If $\bar r$ denotes the class of the conjugate then $\bar r^L=r^R$ and
$\bar r^R=r^L$. 
\item If there is a sector with left vacuum $\mu$ and right vacuum $\nu$
and there is one with left vacuum $\lambda$ and right vacuum $\mu$ then there
is a sector with left vacuum $\lambda$ and right vacuum $\nu$.
\end{itemize}

\section{Regular actions of $C^*$-weak Hopf algebras} \label{s modalg}

\subsection{Module algebras}

Categorically speaking a module algebra $M$ over a WHA or WBA $A$ is a monoid
in the  category of left $A$-modules. More explicitely, 
\begin{enumerate}
\item $M$ is a left $A$-module, the action of $a$ on $m$ is denoted by
$\ma{a}m$, thus
\beanon
\ma{a}(\ma{b}m)&=&\ma{ab}m\,,\qquad a,b\in A,\ m\in M\\
\ma{1}m&=&m\,,\qquad m\in M\,.
\eeanon
\item $M$ is an algebra with unit $1_M$.
\item multiplication of $M$ is an $A$-module map, i.e.,
$$
\ma{a}(mm')\ =\ (\ma{a\1}m)(\ma{a\2}m')\,,\qquad a\in A,\ m,m'\in M\,.
$$
\item the unit of $M$ is an $A$-module map, i.e.,
$$
\ma{a}1_M\ =\ \ma{\PL(a)}1_M\,,\qquad a\in A\,.
$$
\end{enumerate}
If $A$ is a $C^*$-WHA one adds the requirements
\begin{enumerate}
\setcounter{enumi}{4}
\item $M$ is a $C^*$-algebra.
\item $m\mapsto \ma{a}m$ is continuous for all $a\in A$.
\item $(\ma{a}m)^*\ =\ \ma{S(a)^*}m^*$ for $a\in A$, $m\in M$.
\end{enumerate}

\subsection{The invariant subalgebra}

The invariants of an $A$-module algebra $M$ are those elements of 
$M$ that transform the same way under the action of $A$ as the identity $1_M$
does. So we define
\beq
M^A\ :=\ \{n\in M\,|\, \ma{a}n=\ma{\PL(a)}n,\ a\in A\,\}\,.
\eeq
This subspace of $M$ is actually a sub-$C^*$-algebra. It is easy to see that
the Haar measure $h\in A$ provides a conditional expectation onto the
invariant subalgebra, $\ma{h}M=M^A$. 

The structure of general WHA-actions is at least as complicated as the
structure of general actions of finite groups. What we are interested in is
rather the formulation of conditions for a "nice" action (being outer and
Galois, e.g.) which makes the WHA together with its action uniquely determined
by the inclusion of the invariant subalgebra $N=M^A$ in $M$. 

\subsection{The relative commutant}

When $M$ is a module algebra over an ordinary $C^*$-Hopf algebra the
relative commutant $(M^A)'\cap M$ can be arbitrary small. So the "natural
thing" is to consider irreducible inclusions for which $(M^A)'\cap
M=\mathbb{C}1$. There can be situations, however, when an inclusion $N\subset
M$ is reducible and we want to descibe $N$ as an invariant subalgebra. In this
case WHA's are useful for the following reason. For module algebras $M$ over a
$C^*$-WHA there is a non-trivial lower bound for the relative commutant.
Namely, for faithful actions (i.e., $\ma{a}A=0\Rightarrow a=0$) the map
\beq
A^L\ni l\mapsto \ma{l}1_M\in N'\cap M
\eeq
is an injective $C^*$-algebra map. So the "natural thing" is to consider
module algebras for which $(M^A)'\cap M\cong A^L$.

\subsection{The crossed product and Galois actions}

In order to formulate our next condition we need the notion of the crossed
product algebra. The crossed product algebra $M\rtimes A$ of an $A$-module
algebra $M$ with $A$ is equal, as a linear space, to the $A^L$-module tensor
product $M\o_{A^L}A$, where the right $A^L$-module structure of $M$ is
defined by $m\cdot l:=m(\ma{l}1_M)$, and the $^*$-algebra structure on 
$M\o_{A^L}A$ is given by
\bea
(m\rtimes a)(m'\rtimes a')&=&m(\ma{a\1}m')\rtimes a\2 a'\\
(m\rtimes a)^*&=&(\ma{a\1^*}m^*)\rtimes a\2^*\,.
\eea
The crossed product contains $M$ and $A$ as $C^*$-subalgebras in the form
$\{m\rtimes 1\,|\,m\in M\,\}$ and $\{1_M\rtimes a\,|\,a\in A\,\}$,
respectively. In the sequel we identify $M$ and $A$ with these subalgebras. In
this sense we have, for example, the following relation in the crossed product
\beanon
hmh&=&(\ma{h\1}m)h\2 h=(\ma{h\1}m)h\2 S(h\3)h\ =\ (\ma{1\1 h}m)1\2 h\\
   &=&(\ma{h}m)h
\eeanon
proving that the conditional expectation $\ma{h}\under\colon M\to M^A$ is
implemented by the projection $h\in A\subset M\rtimes A$. Therefore the basic
construction $hMh$ is a subalgebra of the crossed product. 

\begin{exa}
If $A$ is a $C^*$-WHA and $\duA$ its dual then $A$ is a $\duA$-module algebra
via the Sweedler arrow, $\ma{\varphi}a=\varphi\la a$, $\varphi\in\duA$, $a\in
A$. The invariant subalgebra is $A^L$. The crossed product $C^*$-algebra
$A\rtimes \duA$ is called the Weyl algebra or Heisenberg double since it is the
algebra generated by $A$ (the "momenta") and $\duA$ (the "coordinates")
satisfying the generalized Weyl commutation relations
\beq
\varphi a\ =\ a\1\bra\varphi\1,a\2\ket\varphi\2\ .
\eeq
It has been shown in \cite{BSz2} that the basic construction for the inclusion
$A^L\subset A$ is precisely the Weyl algebra. 
\end{exa} 

For a Galois action of $A$ on $M$ one requires that the crossed product be
equal to the basic construction for $M^A\subset M$. That is to say
$M\rtimes A$ is generated as an algebra by $M$ and by the Haar element $h$,
being the Jones projection in this case. Thus the above example is a Galois
action.

\begin{exa}
Let $M=A^L$ and let $A$ act on $M$ as the trivial representation:
$\ma{a}l=\PL(al)$. Then $A^L$ is a module algebra in the $C^*$-sense.
The invariant subalgebra is precisely $Z^L$. Now it is clear that the action
of $A$ on $A^L$ is not Galois in general. The crossed product $A^L\rtimes A$
coincides with the $C^*$-algebra $A$. Hence $A^L$ is Galois precisely if
$Z^L\subset A^L\subset A$ is a basic construction. Such special WHA's occur as
symmetries of depth 1 inclusions. 
\end{exa}

\subsection{Regular actions}

The next definition summarizes our requirements on the $A$-module algebra $M$.
\begin{defi}
A module algebra $M$ over the $C^*$-Hopf algebra $A$ is called regular if
\begin{enumerate} \label{def: regular}
\item $M$ is minimal, i.e., $(M^A)'\cap M=A^L$,
\item $M^A\subset M\subset M\rtimes A$ is a basic construction 
in the sense of Jones \cite{G-H-J},
\item the conditional expectation $\ma{h}\under\colon M\to M^A$ is of finite
index in the sense of Watatani \cite{Watatani}.
\end{enumerate}
\end{defi}

For a regular module algebra $M$ one has complete control over the other
relative commutants, too.
\bea
M'\cap (M\rtimes A)&=&A^R\\
(M^A)'\cap (M\rtimes A)&=&A\\
\Center M^A&=&A^L\cap \Center A\\
\Center M&=&A^L\cap A^R\\
\Center (M\rtimes A)&=&A^R\cap\Center A
\eea

The dual WHA $\duA$ acts on the crossed product via the Sweedler arrow on $A$,
i.e.,
\beq
\varphi\blacktriangleright(m\rtimes a)\ :=\ m\rtimes(\varphi\la a)\,,
\quad \varphi\in\duA,\ m\in M,\ a\in A
\eeq
is a left action of $\duA$ on $M\rtimes A$ with the invariant subalgebra being
just $M$. The crossed product algebra $(M\rtimes A)\rtimes\duA$ contains the
algebra $(1_M\rtimes A)\rtimes \duA$ isomorphic to the Weyl algebra. It is
precisely the relative commutant $(M^A)'\cap((M\rtimes A)\rtimes\duA)$.

\section{The reconstruction theorem} \label{s reco}

In this Section we would like to investigate the problem of whether an
inclusion $N\subset M$ of unital $C^*$-algebras is isomorphic to the inclusion
$M^A\subset M$ of the invariant subalgebra with respect to a regular action of
an appropriate $C^*$-WHA $A$. 

If such a WHA-action exists then there must be a conditional expectation
$E\colon$ $M\to N$ of finite index type. Furthermore, $A$ must be
isomorphic, as a $C^*$-algebra, to $N'\cap M_2$ where $M_2$ is the basic
construction for $N\subset M$. 
For obtaining information about the coproduct of $A$ one may look at 
the next member of the Jones tower, $M_3$, because it contains $\duA$ as the 
relative commutant $M\cap M_3$. The derived tower of the Jones tower over
$N\subset M$ is therefore completely known,

\vskip 0.3truecm
\parbox[c]{4.5in}{
\beq
\CD
N@.\subset@.M@.\subset@.M_2@.\subset@.M_3@.\subset@.\dots\\
\cup@.@.\cup@.@.\cup@.@.\cup@.@.\\
N'\cap N@.\ \ \subset\ \ @.N'\cap M@.\ \ \subset\ \ @.N'\cap M_2@.\ \ \subset\
\ @.N'\cap M_3@.\ \ \subset\ \ @.\dots\\
@|@.@|@.@|@.@|@.\\
Z^L@.@.A^L@.@.A@.@.A\rtimes\duA@.@.
\endCD
\eeq
}

\vskip 0.3truecm\noindent
The derived tower is again a Jones tower starting from the second item $A^L$
due to the fact that the Weyl algebra is a Jones extension of $A^L\subset A$
\cite{BSz2}. This means, by definition, that $N\subset M$ is of depth 2. 

Thus the inclusion $N\subset M$ can be the inclusion of the invariant
subalgebra $N=M^A$ w.r.t. a regular $C^*$-WHA action only if it is a finite
index depth 2 inclusion of unital $C^*$-algebras with finite dimensional
centers. The letter condition comes partly from the first derived tower,
saying that $\Center N=Z^L\equiv A^L\cap\Center A$, partly from the second
derived tower, saying that $\Center M=A^L\cap A^R\equiv
\duA^L\cap\Center\duA$. The above conditions are not only neccessary but also
sufficient.
\begin{thm}  \label{thm: reco}
Let $N\subset M$ be an inclusion of unital $C^*$-algebras. Then the
conditions listed below are necessary and sufficient for the existence of
a $C^*$-WHA $A$ and a regular action of $A$ on $M$ such that $M^A=N$.
\begin{enumerate}
\item There exists a conditional expectation $E\colon M\to N$ of finite index
type \cite{Watatani}.
\item $N\subset M$ is of depth 2.
\item $\Center N$ (and/or $\Center M$) is finite dimensional.
\end{enumerate}
The WHA $A$ is uniquely determined by the inclusion if we require that the
restriction of the square of the antipode onto $A^L$ be the identity.
\end{thm}

The proof of this Theorem has not been published yet although it was
implicitely present in \cite{NSzW}. Meanwhile a proof of an important
special case has been appeared in \cite{NV2} where D. Nikshych and L.
Vainerman considered type II$_1$ von Neumann algebra factors. In \cite{NV3}
they went beyond the depth 2 case and have constructed a Galois correspondence
for finite depth inclusions $N\subset M$ of II$_1$ factors. As I will try to
explain below the restriction to factors is not really essential and it just
hides the interesting feature of non-trivial vacuum structure. 

In this section we shall outline a proof
of the above Theorem in the hope of that a detailed proof will be available
in the near future. 

\subsection{A 2-categorical generalization}

An inclusion $N\subset M$ of unital $C^*$-algebras is just a special case of
of a unit preserving $^*$-homomorphism $N\to M$. Unital $C^*$-algebras are the
objects ($=$0-cells) of a $C^*$-2-category the arrows ($=$1-cells) of which
are the unit preserving $^*$-algebra maps and for parallel arrows
$\alpha,\beta\colon N\to M$ the intertwiners ($=$2-cells) from $\alpha$ to
$\beta$ are the elements $t\in M$ satisfying the intertwiner relation
$t\alpha(n)=\beta(n)t$ for all $n\in N$. 
For example the selfintertwiners of the 0-cell $M$, considered as a special
1-cell $\id_M\colon M\to M$, are the elements $c$ of $M$ for which $cm=mc$ for
all $m\in M$, i.e., $\End M=\Center M$ as a $C^*$-algebra.

In general a 2-category $\C$ consists of 0,1, and 2-dimensional cells
$\C^0\subset\C^1\subset\C^2$ and there are two partially defined
associative composition laws for 2-cells $s,t\in\C^2$. 
\begin{itemize}
\item The horizontal composition $s\oh t$ is defined
for 2-cells $s\colon \alpha\to\alpha'$ where $\alpha,\alpha'\colon M\to L$ and
$t\colon\beta\to\beta'$ where $\beta,\beta'\colon N\to M$. 
\item The vertical composition $s\circ t$ is defined for
$s\colon\beta\to\gamma$ and $t\colon\alpha\to\beta$ where
$\alpha,\beta,\gamma\colon N\to M$ are parallel arrows. The result is $s\circ
t\colon \alpha\to\gamma$. 
\end{itemize}
They obey the interchange law
\beq
(s\circ t)\oh(s'\circ t')\ =\ (s\oh s')\circ(t\oh t')
\eeq
for all 2-cells for which the left hand side is defined. The usual convention
is to identify the 0-cells $M$ with the arrows $\id_M\colon M\to M$ 
serving as a (partial) unit for $\oh$ and the 1-cell $\alpha\colon N\to M$
with the identity intertwiner $1_\alpha\colon\alpha\to\alpha$ which serves as
the (partial) unit for $\circ$. So all unis for $\oh$ are units for $\circ$,
too.

The basic example of a 2-category is $\Cat$, the 2-category of small
categories. Its 0-cells are the small categories, the 1-cells are the
functors, and the 2-cells are the natural transformations. 

A 2-category is called a $C^*$-2-category if the sets $\Hom(\alpha,\beta)$ of
2-cells from $\alpha$ to $\beta$ are Banach spaces for all pairs of parallel
1-cells $\alpha,\beta\colon N\to M$, and if the vertical composition makes
$\Hom(\alpha,\alpha)$ into a $C^*$-algebra for each 1-cell $\alpha$.
For a precise definition we refer to \cite{LR}.

The basic example of a $C^*$-2-category is $C^*$-$\Alg$, the category of unital
$C^*$-algebras. Its 0-cells are the small unital $C^*$-algebras, the 1-cells
are the unit preserving $^*$-algebra maps, and the 2-cells are the
intertwiners as it has been defined above.

For our purposes the category $C^*$-$\Alg$ may turn out to be too small in the
following sense. If $N\subset M$ is an inclusion possessing a conditional
expectation $E\colon M\to N$ of finite index, this means that $E$ has a
quasibasis, i.e., a finite set $\{m_i\}$ of elements of $M$ such that
\beq
\sum_{i=1}^n\ m_i\,E(m_i^*m)\ =\ m\,,\qquad \forall m\in M\,.
\eeq
Then one may construct a (2-sided) dual $\bar\iota\colon M\to N$ of the
inclusion map $\iota\colon N\to M$ as a $^*$-algebra map from $M$ to a finite 
amplification $N\o M_n$ of $N$ by the formula
\beq
\bar\iota(m)\ =\ \sum_{i,j}\ E(m_i^*mm_j)\o e_{ij}
\eeq
where $\{e_{ij}\}$ is a set of matrix units of $M_n$.
However, $\bar\iota$ belongs to $C^*$-$\Alg$ only for $n=1$, otherwise it is
an arrow from $M$ to $N$ in the larger $C^*$-2-category $C^*$-$\amp$ of
"amplimorphisms". This category is basically the same as the $C^*$-2-category
of finitely generated projective Hilbert bimodules over unital $C^*$-algebras.
 
We do not need the details about these categories here because we want to
switch to a general $C^*$-2-category $\C$ and lift the conditions (1-2-3)
of Theorem \ref{thm: reco} as assumptions on an arrow $\iota$ of $\C$.
Then we arrive to the following 
\begin{thm} \label{thm main}
Let $\C$ be a $C^*$-2-category and $\iota\colon N\to M$ be an arrow satisfying
the following assumptions.
\begin{enumerate}
\item $\iota$ has a dual (or conjugate) $\bar\iota\colon M\to N$ with
intertwiners $\bar R\colon M\to \iota\oh\bar\iota$ and $ R\colon N\to
\bar\iota\oh\iota$ satisfying the conjugacy equations
\beanon
(\bar R^*\oh\iota)\circ(\iota\oh R)&=&\iota\,,\\
(\bar\iota\oh\bar R^*)\circ(R\oh\bar\iota)&=&\bar\iota\,.
\eeanon
\item $\iota$ is of depth 2, i.e., $\iota\oh\bar\iota\oh\iota$ is a direct
summand of a finite multiple of $\iota$. 
\item $\End N$ is a finite dimensional $C^*$-algebra.
\end{enumerate}
Then $A:=\End\iota\oh\bar\iota$ and $B:=\End\bar\iota\oh\iota$ are finite
dimensional $C^*$-algebras and there exists a non-degenerate bilinear form
$\bra\under,\under\ket\colon A\x B\to\mathbb{C}$ which makes $A$ and $B$ into
$C^*$-WHA's in duality. 
\end{thm}
Under assumptions (1) and (2) assumption (3) is equivalent to assuming $\End
M$ is finite dimensional.

\subsection{Outline of the proof}

\subsubsection{Standard rigidity intertwiners}
Given a conjugate $\bar\iota$ of $\iota$ there is some freedom in choosing the
intertwiners $R$ and $\bar R$. A special class of choices are termed {\em
standard}. They have been defined in case of $C^*$-categories with irreducible
unit in \cite{LR}. It can be generalized easily to $C^*$-categories with
reducible unit $U$ if $U$ decomposes into finitely many irreducibles
\cite{BSz2}. The $C^*$-2-category generalization is straightforward.
Let  \beq \label{eq iota}
\iota\ =\ \bigoplus_a\ \bigoplus_{i=1}^{m_a}\ \iota_a
\eeq
be the decomposition of $\iota$ into pairwise inequivalent irreducibles
$\iota_a\colon N\to M$, each of them with some multiplicity $m_a$. 
Let $\bar\iota_a$ be a conjugate of $\iota_a$ with rigidity intertwiners
$\bar R_a\colon M\to \iota_a\oh\bar\iota_a$ and $ R_a\colon N\to
\bar\iota_a\oh\iota_a$. Then $\bar R_a^*\circ \bar R_a$ is a selfintertwiner 
of $M$ and $ R_a^*\circ R_a$ is a selfintertwiner of $N$. They must be
proportional to a minimal projection of the Abelian algebras $\End M$ and
$\End N$, respectively, since $\iota_a$ is irreducible. Therefore we can
multiply $R_a$ and $\bar R_a$ with numbers, if necessary, to obtain a choice
for which
\beq
\bar R_a^*\circ \bar R_a\ =\ d_a\cdot P_{a^L}\,,\qquad
 R_a^*\circ  R_a\ =\ d_a \cdot P_{a^R}
\eeq
with the same positive number $d_a$ in both equations. Here $P_\mu$ refers to
a minimal (central) projection of $\End M$ or $\End N$ depending on whether the
the $\mu$ is a left or a right "vacuum" of $a$. (Later $\End M$ will become
the $Z^L$ of the WHA $A$ and $\End N$ that of $B$, so this explains the "vacuum
sector" terminology.)

Having been chosen rigidity intertwiners for each $\iota_a$ we are ready to
write down the standard rigidity intertwiners for $\iota$. Let 
\beq
w_{ai}\colon\iota_a\to\iota\,,\quad \bar w_{ai}\colon \bar\iota_a\to\bar\iota
\quad i=1,\dots m_a
\eeq
be intertwiners chosen in such a way that 
\beq
w_{ai}^*\circ w_{bj}=\delta_{ab}\delta_{ij}\ \iota_b\,,\qquad
\sum_a\sum_{i=1}^{m_a}\,w_{ai}\circ w_{ai}^*=\iota\,.
\eeq
I.e., the $w$ intertwiners provide a direct sum diagram for (\ref{eq iota}).
Then the intertwiners
\bea \label{eq standard R}
\bar R&=&\sum_a\sum_i\ (w_{ai}\oh \bar w_{ai})\circ \bar R_a\quad\colon
M\to\iota\oh\bar\iota\\
 R&=&\sum_a\sum_i\ (\bar w_{ai}\oh w_{ai})\circ R_a\quad \colon
N\to\bar\iota\oh\iota
\eea
are rigidity intertwiners for $\iota$. They will be called the standard
rigidity intertwiners. Although they are not unique, depend on the choice of
the direct sum diagram, the maps 
\bea
\Psi_{\iota}\colon \End\iota\to\End M\,,\quad&&\Psi(l):=\bar R^*\circ
(l\oh\bar\iota)\circ \bar R\\
\Phi_{\iota}\colon\End\iota\to\End N\,,\quad&&\Phi(l):=
R^*\circ(\bar\iota\oh l)\circ R
\eea
are uniquely determined faithful positive traces, called the standard traces.
These traces are of finite index type, i.e., have a quasibasis. $\Ind
\Psi_{\iota}=\Ind\Phi_\iota$. Furthermore,
\beq
\tr_M(\Psi_\iota(l))\ =\ \tr_N(\Phi_\iota(l))\,,\quad l\in\End\iota
\eeq
where $tr_M$ is the trace on $\End M$ which takes the value $1$ on each
minimal projector and $\tr_N$ is the analogue trace on $\End N$.

In a similar fashion one defines the standard trace $\Psi_{\bar\iota}\colon
\End\bar\iota\to \End N$, for example, using the standard rigidity intertwiner
$ R$. Also we can construct the standard trace $\Psi_{\iota\oh\bar\iota}$
using the standard rigidity intertwiner
\beq
\bar R_{\iota\oh\bar\iota}\ =\ (\iota\oh R\oh\bar\iota)\circ \bar R\,.
\eeq 
Finally we will need the standard trace $\Psi_{\iota\oh\bar\iota\oh\iota}\colon
\End(\iota\oh\bar\iota\oh\iota)\to \End M$. The following abbreviations will be
used \beq
\Psi_1=\Psi_{\iota}\,,\quad \Psi_{12}:=\Psi_{\iota\oh\bar\iota}\,,\quad
\Psi_{123}:=\Psi_{\iota\oh\bar\iota\oh\iota}\,.
\eeq

\subsubsection{The pairing}
We define the $C^*$-algebras $A:=\End(\iota\oh\bar\iota)$ and
$B:=\End(\bar\iota\oh\iota)$. They are finite dimensional, as all
intertwiner spaces are in a $C^*$-category the monoidal unit of which has
finite dimensional endomorphism algebra. Motivated by the pairing formula of
Theorem 4.11 in \cite{BSz2} we make the following Ansatz
\beq \label{eq: pairing}
\bra a,b\ket\ =\ \tr_M\circ\Psi_{123}(U_{23}U_{12}bz_1az_3)\,,\quad a\in
A,b\in B 
\eeq
where $U_{12}:=\bar R\circ \bar R^*$ and $U_{23}:= R\circ R^*$ are meant
to be embedded into $\End(\iota\oh\bar\iota\oh\iota)$ in the obvious way, just
like $a$ and $b$. The pairing also contains the yet undetermined element
$z\in\End\iota$ in the form of
$z_1:=z\oh\bar\iota\oh\iota$ and $z_3:=\iota\oh\bar\iota\oh z$.

If $z\in\End\iota$ is invertible then the bilinear form (\ref{eq: pairing}) 
is non-degenerate. This can be shown by proving that the Fourier transform
\beq
\F\colon B\to A\,,\quad \F(b):=\Psi_3(U_{23}U_{12}z_1bz_3)
\eeq
is invertible and that the pairing can be written as
\beq
\bra a,b\ket\ =\ \tr_M\circ\Psi_{12}(\F(b)a)\,.
\eeq

By means of the pairing one defines coalgebra structures on $A$ and $B$.
\bea
&\cop_A\colon& A\to A\o A\,,\quad\bra\cop(a),b\o b'\ket=\bra a,bb'\ket\\
&\eps_A\colon& A\to\mathbb{C}\,,\quad \eps_A(a)=\bra a,1_B\ket
\eea
and similar expressions for $B$. Antipodes can be introduced by
\beq
S_A(a):=(a_*)^*\,,\qquad S_B(b):=(b_*)^*
\eeq
where the lower star operation is the transpose of the upper one,
\beq
\bra a_*,b\ket=\overline{\bra a, b^*\ket}\,,\quad
\bra a,b_*\ket=\overline{\bra a^*, b\ket}\,.
\eeq
Explicitely
\bea
\eps_A(a)&=&\tr_M(\bar R^*\circ (z\oh\bar\iota)\circ a\circ 
(z\oh\bar\iota)\circ \bar R)\\
S_A(a)&=&(z\oh z^{-1})\circ \bar a\circ (z^{-1}\oh z)
\eea
where $a\mapsto \bar a$ is the action of the conjugation (left=right duality)
functor w.r.t. the standard conjugacy intertwiners. The explicit form of the
coproduct can be written only in terms of some quasibasis.

The formula for $\eps_A$ shows that $z$ must be Hermitean in order for
$\eps_A$ to be positive. Assuming this the following properties can now be
verified easily: 
\begin{description}
\item[(i)] $(\cop_A\o\id_A)\circ\cop_A\ =\ (\id_A\o\cop_A)\circ\cop_A$
\item[(ii)] $(\eps_A\o\id_A)\circ\cop_A\ =\ \id_A\ =\
            (\id_A\o\eps_A)\circ\cop_A$ 
\item[(iii)] $S_A(aa')\ =\ S_A(a')S_A(a)$
\item[(iv)] $\cop_A\circ S_A\ =\ (S_A\o S_A)\circ\cop_A^{op}$
\item[(v)] $^*\,\circ S_A\circ\,^*\,\circ S_A\ =\ \id_A$
\item[(vi)] $\cop_A(a^*)\ =\ \cop_A(a)^*$
\end{description}
and analogue statements for $B$.

Comparing these properties with the original $C^*$-WHA axioms of \cite{BSz} we
see that what is missing for $A$ with $\cop_A$ to be a $C^*$-WHA, and $B$ with
$\cop_B$ its dual, is the verification of three more axioms:
\begin{description}
\item[(vii)] $\cop_A(a)\cop_A(a')\ =\ \cop_A(aa')$
\item[(viii)] $a\1\o a\2 S_A(a\3)\ =\ 1\1 a\o 1\2$
\item[(ix)] $\eps_A(aa')\ =\ \eps_A(a 1\2)\eps_A(1\1 a')$
\end{description}
The difficult part will be to prove {\bf (vii)} using the depth 2 property and
an appropriate choice of $z$, the remaining ones will hold true without any
further assumptions.

\subsubsection{The depth 2 condition}
If $\iota\oh\bar\iota\oh\iota$ is a direct summand of a finite multiple of
$\iota$ then there exists a finite set of intertwiners $V_i\colon
\iota\to\iota\oh\bar\iota\oh\iota$ such that $\sum_i V_i\circ V_i^*$ is the
identity at $\iota\oh\bar\iota\oh\iota$. Introducing 
\beq
v_i\ :=\ (\iota\oh\bar\iota\oh \bar R^*)\circ (V_i\oh\bar\iota)\ \in\ A
\eeq
we obtain the relation
\beq
\sum_i\ (v_i\oh\iota)\circ U_{23}\circ (v_i^*\oh\iota)\ =\
\iota\oh\bar\iota\oh\iota\,.
\eeq
Hence $\{v_i\}$ is a quasibasis for
$\Psi_2\colon\End(\iota\oh\bar\iota)\to\End\iota$ and
$\End(\iota\oh\bar\iota\oh\iota)$ is the basic construction for
$\End\iota\subset\End(\iota\oh\bar\iota)$. 

\subsubsection{Multiplicativity of the coproduct}
Given a Hermitean invertible $z$ in the pairing of $A$ and $B$ we have the
\begin{pro}
Under the conditions of Theorem \ref{thm main} the following statements are
equivalent: 
\begin{enumerate}
\item $\cop_A$ is multiplicative
\item $\cop_B$ is multiplicative
\item $b\la (aa')\ =\ (b\1\la a)(b\2\la a')$
\item $(\iota\oh b)\circ (a\oh\iota)\ =\ 
((b\1\la a)\oh\iota)\circ(\iota\oh b\2)$
\item if $f\colon\bar\iota\oh\iota\oh\bar\iota\oh\iota\to\bar\iota\oh\iota$
denotes the "fork" $\bar\iota\oh(\bar R^*\circ(z^{-1}\oh\bar\iota))\oh\iota$ 
then 
\beq \label{eq fork rule}
b\circ f\ =\ f\circ (b\1\oh b\2)\,.
\eeq
\end{enumerate}
\end{pro}

Let $\{u_i\}$ be an orthonormal basis of $A$ w.r.t. $\tr_M\circ\Psi_{12}$,
i.e.,
\beq
\tr_M\circ\Psi_{12}(u_i^*u_j)\ =\ \delta_{ij}\,.
\eeq
Then $u^i:=\F^{-1}(u_i^*)\in B$ is the dual basis of $\{u_i\}$ w.r.t.
the pairing. Hence
\beq
\cop_B(b)\ =\ \sum_i\ (u_i\la b)\o u^i\ .
\eeq
\begin{lem}
The fork rule (\ref{eq fork rule}) holds if and only if
$v_i:=u_i^*\circ(z^{-1}\oh\bar\iota)$ is a quasibasis for $\Psi_2\colon
A\to\End\iota$, i.e., if
\beq
\sum_i\ u_i^*z_1^{-2} U_{23}u_i\ =\ 1
\eeq
holds in $\End\iota\oh\bar\iota\oh\iota$.
\end{lem}
The $z's$ that satisfy all conditions we have are precisely the Hermitean
invertible $z\in\End\iota$ for which
\beq
\tr_a z^{-2} \ =\ d_a
\eeq
for all sector $a$ contained in $\iota$. If we add the somewhat ad hoc
requirement that $z$ be central in $\End\iota$ then $z$ is unique up to a sign
in each sector. The $\pm$ signs can be reabsorbed into the freedom of choosing
the standard rigidity intertwiners $R$, $\bar R$ but non-central $z$ cannot be
made central by this method. 

We remark that the $z$ we use here is related to the canonical grouplike
element of the resulting WHA $A$ by the formula $z^2=g'_L$ where the $g'_L$
has been defined in \cite{BSz2}.

Although any $C^*$-WHA can occur as the symmetry of an inclusion $N\subset
M$ satisfying the conditions of Theorem \ref{thm: reco}, we do not need all of
them. Namely, every such inclusion has a uniquely determined $C^*$-WHA with
the additional property $S^2|_{A^L}=\id_{A^L}$. Therefore the freedom of having
non-trivial $S^2|_{A^L}$, or equivalently, non-tracial $\eps|_{A^L}$ in a WHA
is a property which is not utilized in applications to depth 2 inclusions.
Thus the meaning of this degree of freedom of WHA's is still waiting for
explanation.


\begin{thebibliography}{ABC}


\bibitem{BSz} B\"ohm, G., Szlach\'anyi, K., \textit{A coassociative
$C^*$-quantum group with nonintegral dimensions}, Lett. Math.
Phys. \textbf{35} (1996), 437--456 

\bibitem{BNSz} B\"ohm, G., Nill, F., Szlach\'anyi, K. \textit{Weak Hopf
Algebras I:  Integral Theory and the $C^*$-structure},
J. Algebra \textbf{221} (1999), 385--438 

\bibitem{BSz2} B\"ohm, G., Szlach\'anyi, K. \textit{Weak Hopf Algebras II:
Representation theory, dimensions, and the Markov trace}, 
J. Algebra {\bf 233} (2000), 156--212  

\bibitem{EN} Etingof, P., Nikshych, D. \textit{Dynamical quantum groups at
roots of 1}, e-print: \hfill\newline math.QA/0003221

\bibitem{EV} Enock, M., Vallin, J.-M. \textit{Inclusions of von Neumann
algebras and quantum groupoids}, Inst. de Math. de
Jussieu, preprint No.156, 1998

\bibitem{E2} Enock, M. \textit{Inclusions of von Neumann
Algebras and quantum groupoids II}, Inst. de Math. de
Jussieu, preprint No.231, 1999

\bibitem{E3} Enock, M. \textit{Sous-facteurs interm\'ediaires et groupes
quantiques mesur\'es}, J. Operator Theory \textbf{42} (1999), 305--330  
  
\bibitem{G-H-J} Goodman, F.M., de la Harpe, P., Jones, V.F.R.:
\textit{Coxeter Graphs and Towers of Algebras}, Springer 1989

\bibitem{Hay} Hayashi, T. \textit{Quantum group symmetry of partition 
functions of IRF models and its application to Jones' index 
theory}, Commun. Math. Phys. \textbf{157} (1993), 331-345  

\bibitem{LR} Longo, R., Roberts, J. E. \textit{A theory of
dimension}, K-theory \textbf{11} (1997), 103--159 

\bibitem{Lu} Lu, J. H. \textit{Hopf algebroids and quantum groupoids}, 
Int. J. Math. \textbf{7} (1996), 47--70 

\bibitem{NV1} Nikshych, D., Vainerman, L. \textit{Algebraic versions of a
finite-dimensional quantum groupoid},
e-print: math.QA/9808054

\bibitem{NV2} Nikshych, D., Vainerman, L. \textit{A characterization of depth 2
subfactors of $II_1$ factors}, to appear in J. Func. Anal., 
e-print: math.QA/9810028

\bibitem{NV3} Nikshych, D., Vainerman, L. \textit{A Galois correspondence for
II$_1$ factors and quantum groupoids}, e-print: math.QA/0001020

\bibitem{NV4} Nikshych, D., Vainerman, L. \textit{Finite quantum groupoids and
their applications}, e-print: math.QA/0006057 

\bibitem{Nill} Nill, F. \textit{Axioms for Weak Bialgebras},
e-print: math.QA/9805104

\bibitem{NSzW} Nill, F., Szlach\'anyi, K., Wiesbrock, H.-W. \textit{Weak Hopf
algebras and reducible Jones inclusions of depth 2},
e-print: math.QA/9806130

\bibitem{Szlach} K. Szlach\'anyi, \textit{Finite quantum groupoids and
inclusions of finite type}, Fields Inst. Commun. Series. to appear, e-print: 
math.QA/0011036

\bibitem{Val1} Vallin, J.-M. \textit{Bimodules de Hopf et poids op\'eratoriels
de Haar}, J. Operator Theory \textbf{35} (1996), 39--65 

\bibitem{Watatani} Watatani, Y. \textit{Index for $C^*$-subalgebras},
Memoirs of the Amer. Math. Soc., No. 424, 1990

\bibitem{Yam} Yamanouchi, T. \textit{Duality for generalized Kac algebras
and a characterization of finite groupoid algebras},
J. Algebra \textbf{163} (1994), 9-50 

\end{thebibliography}
\end{document}